\documentclass[12pt]{article}
\usepackage{graphicx}
\usepackage{amsfonts}
\usepackage{latexsym}
\usepackage{amsmath}

\makeatletter
\@addtoreset{figure}{section}
\def\thefigure{\thesection.\@arabic\c@figure}
\@addtoreset{table}{section}
\def\thetable{\thesection.\@arabic\c@table}

\def\@sect#1#2#3#4#5#6[#7]#8{\ifnum #2>\c@secnumdepth
     \def\@svsec{}\else
     \refstepcounter{#1}\edef\@svsec{\csname the#1\endcsname.\hskip .75em
}\fi
     \@tempskipa #5\relax
      \ifdim \@tempskipa>\z@
        \begingroup #6\relax
          \@hangfrom{\hskip #3\relax\@svsec}{\interlinepenalty \@M #8\par}%
        \endgroup
       \csname #1mark\endcsname{#7}\addcontentsline
         {toc}{#1}{\ifnum #2>\c@secnumdepth \else
                      \protect\numberline{\csname the#1\endcsname}\fi
                    #7}\else
        \def\@svsechd{#6\hskip #3\@svsec #8\csname #1mark\endcsname
                      {#7}\addcontentsline
                           {toc}{#1}{\ifnum #2>\c@secnumdepth \else
                             \protect\numberline{\csname the#1\endcsname}\fi
                       #7}}\fi
     \@xsect{#5}}
\def\@begintheorem#1#2{\it \trivlist \item[\hskip \labelsep{\bf #1\ #2.}]}
\def\section{\@startsection {section}{1}{\z@}{-3.5ex plus -1ex minus
 -.2ex}{2.3ex plus .2ex}{\normalsize\bf}}

\begin{document}

\title{Minimum Perimeter Rectangles That Enclose\\
Congruent Non-Overlapping Circles}
\date{} 
\maketitle

\begin{center}
\author{Boris D. Lubachevsky  \ \ \ \ \ \ \ \ \ \ \ \ \ \ \ \ \ \ \ \ \ \ Ronald L. Graham\\
{\em lubachevsky@netscape.net \ \ \ \ \ \ \ \ \ \ \ \ \ \ \ \ \ \ graham@ucsd.edu}\\
Bell Laboratories \ \ \ \ \ \ \ \ \ \ \ \ \ \ \ \ \ \ \ \ \ \ University of California \\
600 Mountain Avenue  \ \ \ \ \ \ \ \ \ \ \ \ \ \ \ \ \ \ \ \ \ \ \ \ \ \ \ \ \ at San Diego \\
Murray Hill, New Jersey \ \ \ \ \ \ \ \ \ \ \ \ \ \ \ \ \ \ La Jolla, California }
\end{center}

\setlength{\baselineskip}{0.995\baselineskip}
\normalsize
\vspace{0.5\baselineskip}
\vspace{1.5\baselineskip}

\begin{abstract}
We use computational experiments 
to find the rectangles of minimum perimeter into which a given number
$n$ of non-overlapping congruent circles can be packed. No
assumption is made on the shape of the rectangles. 
In many of the packings found, 
the circles form the usual regular square-grid or hexagonal
patterns or their hybrids. However,
for most values of $n$ in the tested
range $n \le 5000$, e.g., for
$n = 7, 13, 17, 21, 22, 26, 31, 37, 38, 41, 43...,4997, 4998, 4999, 5000$,
we prove that the optimum cannot possibly be achieved
by such regular arrangements.
Usually, the irregularities in 
the best packings found for such $n$ are
small, localized modifications to regular patterns;
those irregularities are usually easy to predict.
Yet
for some such irregular $n$,
the best packings found show substantial, extended irregularities
which we did not anticipate.
In the range we explored carefully, 
the optimal packings were substantially irregular 
only for $n$ of the form
$n = k(k+1)+1$, $k = 3, 4, 5, 6, 7$, 
i.e. for $n = 13, 21, 31, 43$, and 57.
Also, we prove that the height-to-width ratio
of rectangles of minimum perimeter containing packings of $n$ congruent
circles tends to 1 as $n \rightarrow \infty$.

{\bf Key words}: disk packings, rectangle, container design,
hexagonal, square grid

{\bf AMS subject classification:} primary 52C15, secondary 05B40, 90C59

\end{abstract}
\section{Introduction}\label{sec:intro}
\hspace*{\parindent} 
Consider the task of finding
the rectangular region of least perimeter
that encloses a given number $n$ of 
circular disks of equal diameter.
The circles must not
overlap with each other or extend outside the rectangle.
The 
aspect ratio of the rectangle,
i.e. the ratio of its height to width,
is variable and 
subject to the perimeter-minimizing choice
as well as the positions of the circles
inside the rectangle.

Dense packings of circles in rectangles of a fixed shape,
in particular, in squares,
have been the subject of
many investigations \cite{GL2}, \cite{NO1}, \cite{NO2}, \cite{NO3};
a comprehensive survey is given in \cite{SMC},
see also \cite{Specht}.
When the aspect ratio of the enclosing rectangle is fixed,
the densest packing minimizes 
both the area and the perimeter of the rectangle.
However, by allowing a variable aspect ratio, 
the two optima may differ for the same number of circles $n$.
In 1970, 
Ruda relaxed the restriction of the fixed aspect
ratio while trying to minimize the area, see
\cite{Ruda}.
He found the minimum-area packings of $n$ congruent circles
in the variably shaped rectangles for $n \le 8$ and
conjectured the minima for $9 \le n \le 12$.
In \cite{LG}, we extended Ruda's conjectures to
$n \le 5000$. 

In this paper we switch our attention to 
minimizing the rectangle perimeter,
while keeping the aspect ratio of the rectangle variable.
We report the results of essentially 
the same computational procedure for finding the minimum-perimeter
packings as the procedure used in \cite{LG} for finding 
the minimum-area packings.
Even though
the optima themselves are usually different,
the structures of the set of optimal
packings turn out to be
very similar 
for the two minimization tasks.
In either case for many $n$,
the optimum pattern
is regular, i.e. it is a square-grid pattern,
or a hexagonal pattern, or a hybrid of these two patterns.
One difference is that the occurrence of 
the non-regular patterns is more frequent in the minimum-perimeter
case than in the minimum-area case:
the smallest non-regular $n$ for the minimum-perimeter criterion is $n = 7$ 
while that for the minimum-area criterion is $n = 49$ in \cite{LG};
for almost all $n$ that are close to $n = 5000$
the minimum-perimeter packings are not regular
while
for the majority of $n$ everywhere in the range $1 \le n \le 5000$,
the regular packings 
still supply the minimum-area rectangles.

It appeared in \cite{LG},
that
if the minimum-area rectangular packing 
for a particular $n \le 5000$ is not regular,
then
the minimum could be obtained by a small and 
easy-to-predict, 
localized modification to a regular pattern.
In the case of the minimum perimeter,
modifications of the same type
apparently supply the optima for most
numbers of circles $n$ in the studied range 
$n \le 5000$ for which the minimum-perimeter rectangular
packing happens not to be regular.
But not for all such numbers!
For certain exceptional numbers of circles $n$,
the packing pattern with the smallest perimeter we found is
complex and/or requires extended modifications to a regular
pattern.\footnote{
We now realize that similar exceptions might also exist in the case of minimizing
the area.
See the footnote
in Section~\ref{sec:largn} for an example.
}
A further surprise was that
these exceptionally irregular packings, 
their patterns being complex and unpredictable, 
seem to occur quite predictably. In particular,
in the range $n \le 62$ which we explored carefully,
our experiments detected such irregular packings
only for $n$ of the form
$n=k(k+1)+1$, i.e.  
for $n = 13, 21, 31, 43$, and $57$. 
The best packings
for larger terms $n$ of this sequence, 
i.e. for $n = 73$, 91,..., are probably similarly irregular,
although we could not test that as thoroughly 
as for the smaller terms
because the computational resources
needed for such a testing grow very rapidly with $n$.

Most of our findings are unproven conjectures, 
the outcomes of computer experiments. 
We do not know why the exceptionally irregular values of $n$
appear along the sequence $n=k(k+1)+1$
and only speculate
by suggesting a possible reason.
A few facts which we can prove are 
explicitly stated as being provable.

\section{Computational method}\label{sec:cmethod}
\hspace*{\parindent} 
To obtain the minimum-perimeter packing conjectures
we use a variant of
the computational technique employed in \cite{LG}
for generating the minimum-area packing conjectures.
The technique consists of two independent algorithms:
the restricted search algorithm and 
the ``compactor'' simulation algorithm.
We now review these procedures.

The restricted search algorithm operates on the
assumption
that the desired minimum is achieved on a set of 
configurations which is much smaller than the set of all possible
configurations. 
The set is restricted to include only
hexagonal patterns,   square-grid  patterns, their hybrids,
and the patterns obtained by removing some circles from
these patterns. 
For a given number of $n > 0$ circles, 
a configuration in the restricted set $R_n$
is defined by 6 integers:
\\
$~~~~w$, the number of circles in the longest row, 
\\
$~~~~h$, the number of rows arranged in a hexagonal alternating pattern, 
\\
$~~~~h_{-}$, those among the $h$ hexagonally arranged rows 
that consist of $w-1$ circles each; 
the rest $h - h_-$ rows consist of $w$ circles each,
\\
$~~~~s$, the number of rows, in addition to $h$ rows, 
that are  stacked in the square-grid pattern, 
\\
$~~~~s_-$, those among the $s$ square-grid rows that consist of $w-1$
circles each; the remaining $s - s_-$ rows consist of $w$ circles each,
\\
$~~~~v$, the number of ``mono-vacancies'' or holes.
\\
The numbers must be non-negative and must
satisfy the following additional restrictions:
\\
$w > 0$, $h + s > 0$, $s_- \le s$, $s_- < s + h$, $h \ne 1$,
$v \le  \min \{ w,h + s \} - 1$,
\\
if $h$ is even, $h = 2k$, 
then $h_-$ can take on only two values $h_- = 0$ or $h_- = k$,
\\
if $h$ is odd, $h = 2k+1$, 
then $h_-$ can take on only three values $h_- = 0$ or $h_- = k$ or $h = k+1$.
\\
Finally, the total number of circles must equal $n$,
\begin{equation}
\label{totaln}
w (h + s) - h_- - s_- - v = n .
\end{equation}

\begin{figure}
\centering
\includegraphics*[width=3.45in]{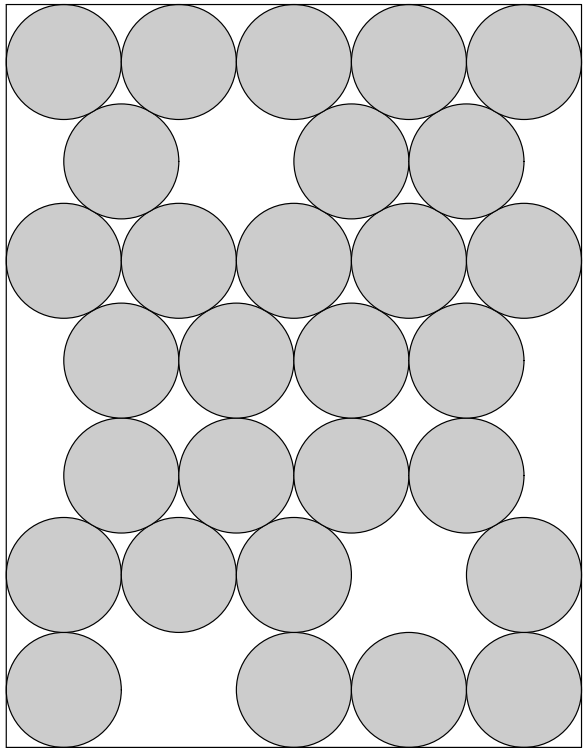}
\caption{A configuration in the restricted set $R_n$ 
for $n=29$ circles;
here $w = 5$, $h = 5$, $h_- = 2$, $s = 2$, $s_{-} = 1$, and $v = 3$,
so that $29 = w(h + s) - h_-  - s_{-} - v$.
}
\label{fig:class}
\end{figure}

A ``general case'' example is shown in Figure \ref{fig:class}. 
As further examples,
we now identify the 6-tuples ($w$, $h$, $h_-$, $s$, $s_-$, $v$)
for some configurations presented in the following sections.
In each example that follows, the unmentioned parameters 
of the tuple are equal zero.
   
In Figure \ref{fig:same}: 
\\
configuration "1 circle" has $w = s = 1$;
\\
configuration ``2 circles'' has $w = 2$, $s = 1$; 
it can also be identified as having $w = 1$, $s = 2$;
\\
configuration ``4 circles'' has $w = s = 2$;
\\
configuration ``6 circles'' has $w = 3$, $s = 2$;
it can also be identified as having $w = 2$, $s = 3$;
\\
configuration ``9 circles'' has $w = s = 3$;
\\
configuration ``11 circles'' has $w = 4$, $h = 3$, $h_- = 1$;
\\
configuration ``12 circles'' has $w = 4$, $s = 3$;
it can also be identified as having $w = 3$, $s = 4$;
\\
configuration ``15 circles'' has $w = 4$, $h = 3$, $h_- = 1$, $s = 1$;
\\
configuration ``19 circles'' has $w = 5$, $h = 3$, $h_- = 1$, $s = 1$.
  
In Figure \ref{fig:7}: 
\\
configuration $a$ has $w = h = 3$, $h_- = 2$;
\\
configuration $b$ has $w = 2$, $h = 3$, $h_- = 1$, $s = 1$;
\\
configuration $c$ has $w = h = 3$, $h_- = 1$, $v = 1$.
  
In Figure \ref{fig:17a26}: 
\\
configuration $a$ has $w = 4$, $h = 5$, $h_- = 2$, $v = 1$;
\\
configuration $c$ has $w = 5$, $h = 6$, $h_- = 3$, $v = 1$.

In Figure \ref{fig:200}: 
\\
configuration $a$ has $w = 13$, $h = 16$, $h_- = 8$;
\\
configuration $b$ has $w = 29$, $h = 7$, $h_- = 3$.
 
Given the values of $w$, $h$, $h_-$, and $s$, 
and that of the common radius of the circles $r$,
the height $H$, width $W$, and perimeter $P$ 
of the enclosing rectangle can be found from
\begin{equation}
\label{height}
\ \ \ H/r  = 2 s + \left\{  \begin{array}{ll}
               2 + ( h - 1 ) \sqrt 3 &   \mbox{if $h > 0$} \\
               0                     &   \mbox{if $h = 0$}
              \end{array}
        \right .
\end{equation}
\begin{equation}
\label{width}
W/r  = 2 w + \left\{  \begin{array}{ll}
               1    &   \mbox{if $h > 0$ and $h_{-} = 0$} \\
               0    &   \mbox{if $h = 0$ \ or \ $h_{-} > 0$}
              \end{array}
        \right .
\end{equation}
\begin{equation}
\label{perim}
P = 2 ( W + H )
\end{equation}
Note that this ratio $P/r$ is a number of the form $ x + y \sqrt 3 $,
where the non-negative integers
$x$ and $y$ are obvious functions
of $w,h,h_-$ and $s$.

Sometimes several different patterns correspond
to a given tuple $w$, $h$, $h_-$, $s$, $s_-$, $v$;
those may differ in the ways the $s$ or $s_-$ rows are attached or
the $v$ holes are selected.
Since the shape and the perimeter of the enclosing rectangle
do not change among these variations, 
the minimization procedure treats them all as the same packing.

An important fact is that for each $n > 0$,
there are only a finite number of 6-tuples 
$w$, $h$, $h_-$, $s$, $s_-$, $v$,
that satisfy the restrictions above. 
For a given value of $n$,
our procedure lists all such 6-tuples,
and for each of them computes $P/r$,
and selects the configurations that correspond to the minimum
value of $P/r$.
By presenting the values $P/r$ in the form $x+y\sqrt{3}$
with integers $x$ and $y$,
only comparisons among integers are involved
and the selection of the minimum is exact.

The reader should be reminded 
that we do not claim that our restricted
search procedure produces
the global optimum for packing $n$ circles in a rectangle. 
In fact, we include 
some configurations which are clearly
non-optimal,
for example, those with $v > 0$.
The usefulness of the given definition of the sets $R_n$
should become apparent when we compare below the restricted
search algorithm outcomes with those of the
``compactor'' simulation algorithm. 

The ``compactor'' simulation works as follow
(see also \cite{PAS}).
It begins by generating a random starting configuration 
with $n$ circles lying inside a (large) rectangle
without circle-circle overlaps. 
The starting configuration is feasible but is usually rather sparse.
Then the computer imitates a ``compactor'' with
each side of the rectangle pressing against the circles,
so that the circles are being forced
towards each other until they ``jam.''
Possible circle-circle or circle-boundary conflicts
are resolved using a simulation of hard collisions so
that no overlaps or boundary-penetrating 
circles occur during the process.

The simulation for a particular $n$ is repeated many times,
with different starting circle configurations.
If the final perimeter value in a run is smaller than
the record achieved thus far, it replaces the current record.
Eventually in this process, the record stops improving up to
the level of accuracy allowed by the double precision accuracy
of the computer. The resulting packing now becomes a candidate
for the optimal packing for this value of $n$.

The main advantage of the ``compactor'' simulation vs.
the restricted search is that in the simulation
no assumption is made about the resulting packing pattern.
The circles are ``free to choose'' any final configuration
as long as it is ``jammed.'' 
This advantage comes at a price: the simulation time needed
in multiple attempts
to achieve a good candidate packing for a particular $n$
is typically several orders of magnitude longer than the time needed
on the same computer to deliver
the minimum in set $R_n$ by the restricted search procedure.
For example, it may take a fraction of a second to 
find the minimum perimeter
packing of 15 circles by the search in $R_{15}$ and 
it may take days
with thousands of attempts
to produce the same answer by simulation.

\section{Results: regular and semi-regular optimal packings}\label{sec:regasem}
\hspace*{\parindent} 
Table \ref{tab:1t62} lists the packings of $n$ equal circles
in rectangles of the smallest found perimeter 
for each $n$ in the range
$1 \le n \le 62$, except $n = 13, 21, 31, 43$, and $57$.
A somewhat arbitrary bound $n = 62$ was set so that
the ``compactor'' simulation was performed
for each $n \le 62$ but only
for a few isolated values $n > 62$
since the simulation slows down significantly for larger $n$.
On the other hand, 
the minimum perimeter packings in $R_n$ were produced
by the restricted search procedure for each $n \le 5000$.

\begin{table}
\begin{center}
\fbox{
\begin{tabular}{r|r|r|r|r|r||r|r|r|r|r|r||r|r|r|r|r|r} 
$n$&$w$&$h$&$h_{-}$&$s$&$\delta$&
$n$&$w$&$h$&$h_{-}$&$s$&$\delta$&
$n$&$w$&$h$&$h_{-}$&$s$&$\delta$ \\ \hline
   1  & 1 & 0 & 0 & 1 & 0         & *22 & 5 & 5 & 2 & 0 &$\delta_2$
                                    & *41 & 6 & 7 & 0 & 0 &$\delta_2$ \\
   2  & 2 & 0 & 0 & 1 & 0         &  23 & 5 & 5 & 2 & 0 & 0
                                    &  42 & 6 & 7 & 0 & 0 & 0         \\
   3  & 2 & 2 & 1 & 0 & 0         &  24 & 5 & 3 & 1 & 2 & 0
                                    &  44 & 6 & 8 & 4 & 0 & 0         \\
   4  & 2 & 0 & 0 & 2 & 0         &  25 & 5 & 5 & 0 & 0 & 0
                                    & *45 & 7 & 7 & 3 & 0 &$\delta_3$ \\
   5  & 2 & 3 & 1 & 0 & 0         & *26 & 5 & 6 & 3 & 0 &$\delta_2$
                                    &  46 & 7 & 7 & 3 & 0 & 0         \\
   6  & 3 & 0 & 0 & 2 & 0         &  27 & 5 & 6 & 3 & 0 & 0
                                    &  47 & 7 & 5 & 2 & 2 & 0         \\
  *7  & 3 & 3 & 1 & 0 &$\delta_1$ &  28 & 5 & 5 & 2 & 1 & 0
                                    &  48 & 6 & 8 & 0 & 0 & 0         \\
   8  & 3 & 3 & 1 & 0 & 0         &     & 6 & 5 & 2 & 0 & 0
                                    &  49 & 7 & 7 & 0 & 0 & 0         \\
   9  & 3 & 0 & 0 & 3 & 0         &  29 & 5 & 3 & 1 & 3 & 0
                                    &  50 & 6 & 9 & 4 & 0 & 0         \\
  10  & 3 & 4 & 2 & 0 & 0         &     & 6 & 3 & 1 & 2 & 0
                                    & *51 & 7 & 8 & 4 & 0 &$\delta_3$ \\
  11  & 3 & 3 & 1 & 1 & 0         &  30 & 5 & 6 & 0 & 0 & 0
                                    &  52 & 7 & 8 & 4 & 0 & 0         \\
      & 4 & 3 & 1 & 0 & 0         &  32 & 5 & 7 & 3 & 0 & 0
                                    &  53 & 7 & 7 & 3 & 1 & 0         \\
  12  & 4 & 0 & 0 & 3 & 0         &  33 & 6 & 6 & 3 & 0 & 0
                                    &     & 8 & 7 & 3 & 0 & 0         \\
  14  & 4 & 4 & 2 & 0 & 0         &  34 & 6 & 5 & 2 & 1 & 0
                                    &  54 & 6 & 9 & 0 & 0 & 0         \\
  15  & 4 & 3 & 1 & 1 & 0         &  35 & 5 & 7 & 0 & 0 & 0
                                    & *55 & 7 & 8 & 0 & 0 &$\delta_3$ \\
  16  & 4 & 0 & 0 & 4 & 0         &  36 & 6 & 6 & 0 & 0 & 0
                                    &  56 & 7 & 8 & 0 & 0 & 0         \\
 *17  & 4 & 5 & 2 & 0 &$\delta_2$ &**37 & 6 & 7 & 3 & 0 &$\delta_1$
                                    & *58 & 7 & 9 & 4 & 0 &$\delta_4$ \\
  18  & 4 & 5 & 2 & 0 & 0         & *38 & 6 & 7 & 3 & 0 &$\delta_3$
                                    &  59 & 7 & 9 & 4 & 0 & 0         \\
  19  & 4 & 3 & 1 & 2 & 0         &  39 & 6 & 7 & 3 & 0 & 0
                                    &  60 & 8 & 8 & 4 & 0 & 0         \\
      & 5 & 3 & 1 & 1 & 0         &  40 & 6 & 5 & 2 & 2 & 0
                                    &  61 & 8 & 7 & 3 & 1 & 0         \\
  20  & 4 & 5 & 0 & 0 & 0         &     & 7 & 5 & 2 & 1 & 0
                                    & *62 & 7 & 9 & 0 & 0 &$\delta_3$ \\
\end{tabular}
}
\caption{Packings of $n$ circles in rectangles 
of the smallest found perimeter
for all $n$ in the range $1 \le n \le 62$, 
except $n = 13, 21, 31, 43$, and 57.
The packing patterns are described with parameters
$w$, $h$, $h_-$, $s$, and $\delta_i$ and with star markings
as explained in the text
}
\label{tab:1t62}
\end{center}
\end{table}

All packings presented 
in Table \ref{tab:1t62} 
can be split into two sets.
The first set consists 
of either perfectly hexagonal packings
or perfectly square-grid packings 
or their hybrids. 
We will call these {\em regular} packings.
A regular packing of $n$ circles is characterized
in Table \ref{tab:1t62} 
by the parameters $n$, $w$, $h$, $h_-$, and $s$,
as defined
in Section~\ref{sec:cmethod}.
Note that the 
parameters $s_-$ and $v$ which are also defined
in Section~\ref{sec:cmethod}
are not present
in Table \ref{tab:1t62}.
These $s_-$ and $v$ equal 0 for each regular packing.

For example, 
Table \ref{tab:1t62} lists
two conjectured minimum-perimeter packings of $n=11$ circles:
one with 
$w=3$, $h=3$, $h_-=1$, and $s=1$,
and the other with 
$w=4$, $h=3$, $h_-=1$, and $s=0$.
The latter is perfectly hexagonal as seen
in Figure \ref{fig:same} (configuration ``11 circles'').
The former is a hybrid of hexagonal and square-grid packings.
A similar hybrid is also 
the only conjectured minimum-perimeter packing of $n=15$ circles
which is listed in
Table \ref{tab:1t62}.
It has parameters
$w=4$, $h=3$, $h_-=1$, and $s=1$
and it is shown 
in Figure \ref{fig:same} (configuration ``15 circles'').

Given the parameters of a regular packing, 
one easily determines the shape and the perimeter of the
rectangle that encloses the circles
using formulas \eqref{height}, \eqref{width}, and \eqref{perim}.
The perimeter for each conjectured minimum-perimeter packings of 11 circles is 
$20 + 4 \sqrt{3}$
and for that of 15 circles
$24 + 4 \sqrt{3}$
in units equal to the common circle radius.

It will be convenient to define the shape of a rectangle
by the ratio $L/S$ where
\begin{equation}
\label{los}
L = \max \{ H, W \},  
S = \min \{ H, W \},  
\end{equation}
so that 
\begin{equation}
\label{argt1}
L/S \ge 1 .
\end{equation}
In these three examples, we have, respectively 
\\
$L/S = (4 + 2 \sqrt{3})/6 
= 1.2440169..$, 
for the first packing of 11 circles
in Table~\ref{tab:1t62};
\\
$L/S = 8/(2 + 2 \sqrt{3}) 
= 1.4541016..$, 
for the second packing of 11 circles
in Table~\ref{tab:1t62};
\\
$L/S = 8/(4 + 2 \sqrt{3}) 
= 1.0717968..$, 
for the packing of 15 circles.

Similar simple calculations can be done for all 
the other regular packings.

All entries $n$ that correspond to regular packings
are not marked by stars in the table.
Now consider the other set of the packings,
those with entries $n$ that are marked by stars in the table.
The smallest example is $n=7$.
\begin{figure}
\centering
\includegraphics*[width=5.1in]{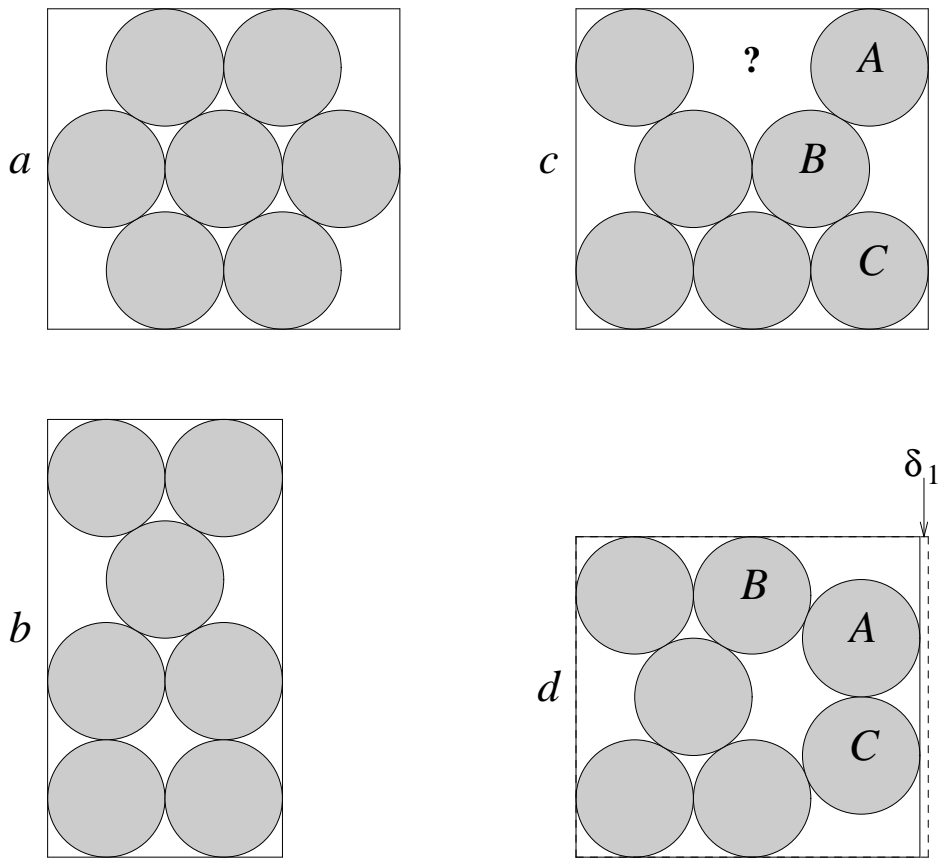}
\caption{
Packings of 7 equal circles in the smallest perimeter rectangle:
$a$) among the hexagonal configurations,
$b$) among the hybrid configurations, 
$c$) among the regular-with-holes configurations,
$d$) among all configurations tested.
The perimeters of the rectangles of the packings $a$, $b$, and $c$ are the same
and are larger than the perimeter of the packing $d$.
}
\label{fig:7}
\end{figure}
The smallest perimeter rectangle that we could find 
(using the simulated ``compactor'' procedure)
is that for configuration $d$
in Figure~\ref{fig:7}.
Because
this configuration does not fit
the description 
in Section~\ref{sec:cmethod}
of a possible pattern in $R_7$,
the restricted search procedure cannot find the packing $d$.
Instead, the best in $R_7$ happen to be
the configurations $a$, $b$, and $c$.
Having among them different aspect ratios 
of the enclosing rectangle,
the configurations $a$, $b$, and $c$ are of the same perimeter
$P/r = 16 + 4 \sqrt{3} = 22.928293..$.
Out of these three, the $c$ requires the smallest
number of circles to be moved 
and the smallest readjustment of the boundary
to obtain the $d$.
We move circles labeled 
in Figure~\ref{fig:7} 
as $A$, $B$, and $C$ 
to turn the $c$ into the $d$.

For the entry $n = 7$, 
Table~\ref{tab:1t62} lists parameters 
$w = 3$, $h = 3$, and $h_- = 1$ and 
those are of the configuration $c$.
Also the entry is marked with one star 
which represents one mono-vacancy
in the pattern $c$.
The implied convention here is that this entry represents
the configuration $d$ because $d$, while
not being describable in the terms of the table,
can be
obtained by a simple and standard transformation
from $c$.
(Note that the configuration $c$ can be defined in several 
ways depending on the position of the hole; the resulting
different configurations are not distinguished
by the restricted search procedure or in the table.)
As $c$ is turned into $d$, the width of the rectangle decreases
by the value $\delta = \delta_1$ where
\begin{equation}
\label{delta1}
\delta_1 = 2 - \sqrt{2 \sqrt {3}}.
\end{equation}
The smallest found perimeter for the 7 circles thus becomes
\begin{equation}
\label{odd}
P^{opt}/r = P/r - 2 \delta
\end{equation}
which is  22.650743.. here.

The $\delta$ will be called the {\em improvement} parameter.
The equality $\delta = 0$ together with
the absence of star markings distinguishes
a regular packing entry in
Table \ref{tab:1t62}.
On the other hand, the entries
with $\delta = \delta_i > 0,~i = 1, 2, 3$ or 4, in
Table \ref{tab:1t62} 
correspond to packings that are not regular.
The number of stars that marks the $n$ for such an entry equals
the number of mono-vacancies $v$
in the packing according to the definition of class $R_n$
in Section~\ref{sec:cmethod}.

\begin{figure}
\centering
\includegraphics*[width=6.0in]{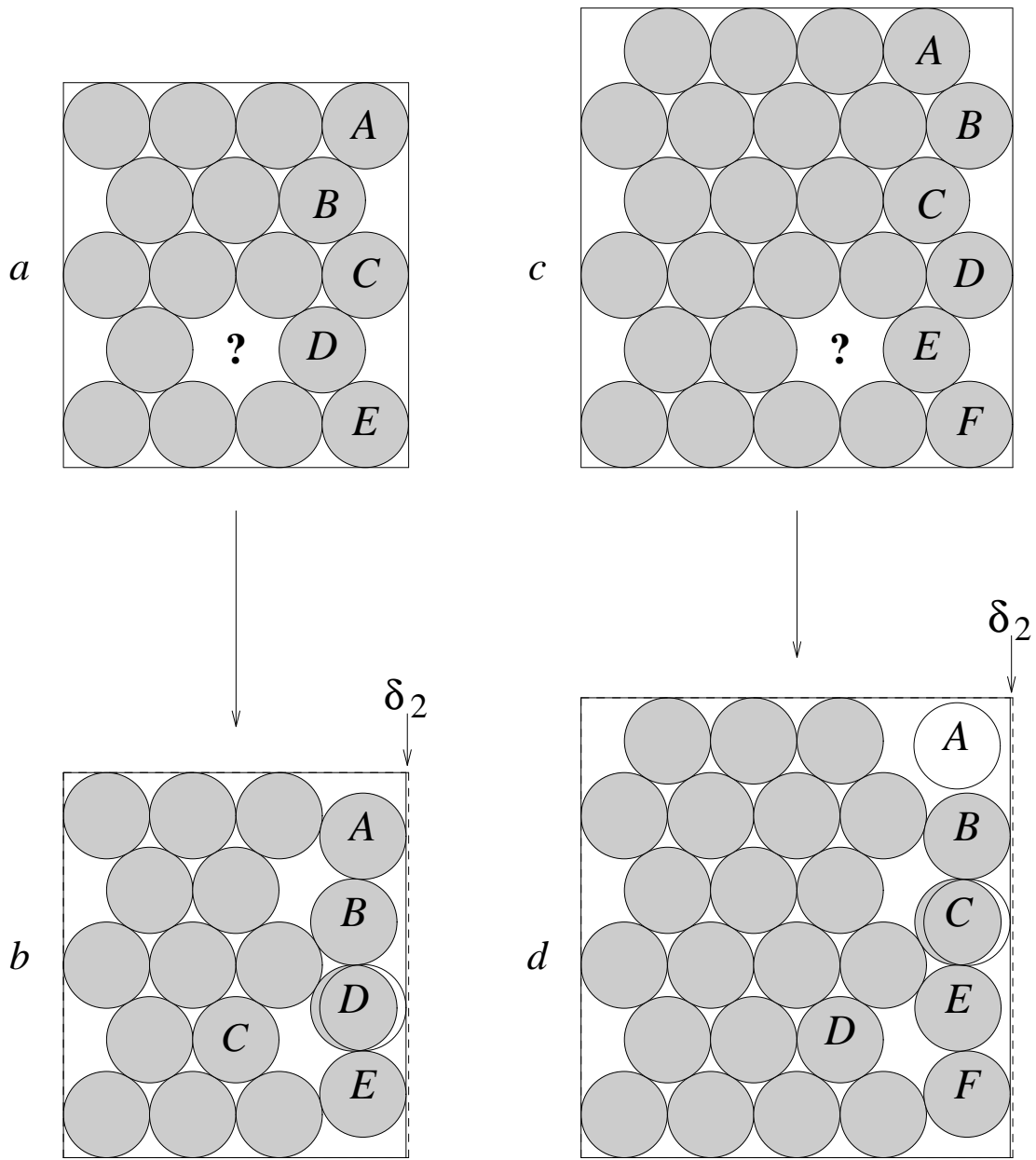}
\caption{Packings of 17 ($a$ and $b$) and 26 ($c$ and $d$)
equal circles in rectangles. Packings $a$ and $c$ are the best
in $R_{17}$ and $R_{26}$, respectively.
Packings $b$ and $d$ are the best we could find
for their number of circles.
Alternative equivalent positions of circles $D$ in packing $b$
and $C$ in packing $d$ are shown 
}
\label{fig:17a26}
\end{figure}
\begin{figure}
\centering
\includegraphics*[width=6.0in]{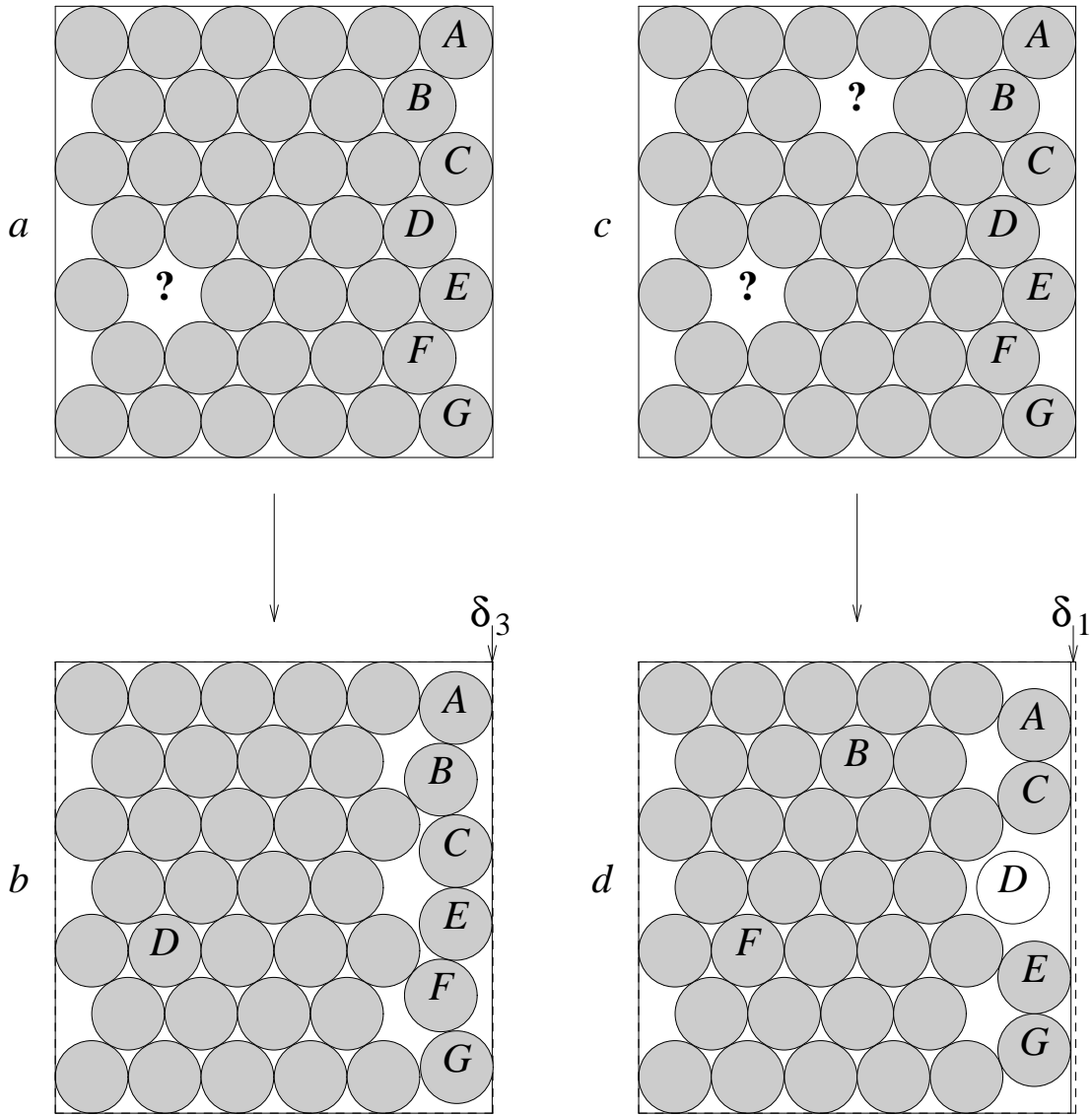}
\caption{Packings of 38 ($a$ and $b$) and 37 ($c$ and $d$)
equal circles in rectangles. Packings $a$ and $c$ are the best
in $R_{38}$ and $R_{37}$, respectively.
Packings $b$ and $d$ are the best we could find
for their number of circles
}
\label{fig:37a38}
\end{figure}
\begin{figure}
\centering
\includegraphics*[width=6.0in]{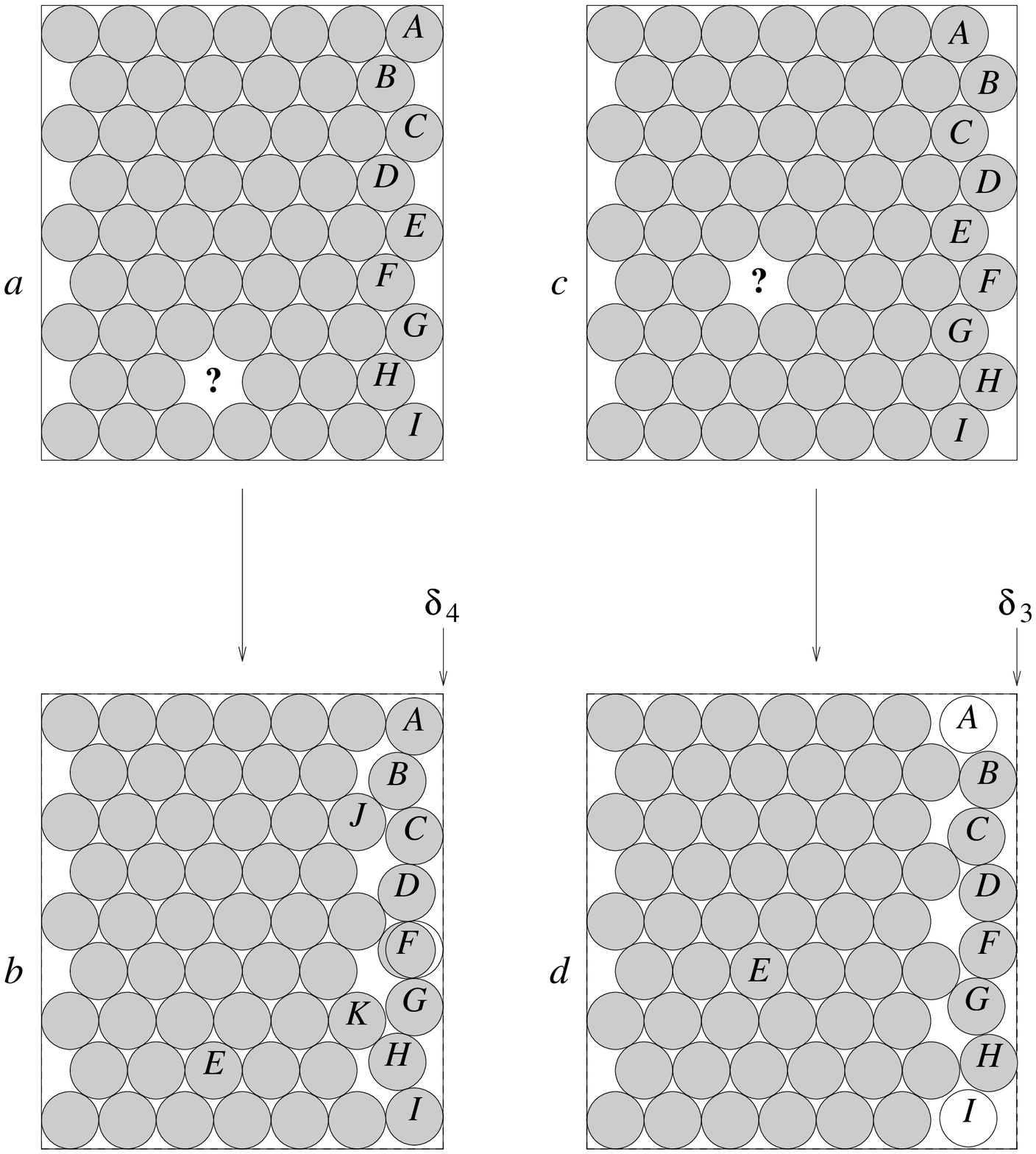}
\caption{Packings of 58 ($a$ and $b$) and 62 ($c$ and $d$)
equal circles in rectangles. Packings $a$ and $c$ are the best
in $R_{58}$ and $R_{62}$, respectively.
Packings $b$ and $d$ are the best we could find
for their number of circles.
A positive gap exists between circles $C$ and $J$
and also between circles $G$ and $K$ in packing $b$.
An alternative equivalent position of circle $F$ is shown
in packing $b$
}
\label{fig:58a62}
\end{figure}
Figures~\ref{fig:17a26}, \ref{fig:37a38}, and \ref{fig:58a62}
show six other
non-regular packings. 
Those are labeled 
$b$ and $d$ in each figure.
Their regular-with-holes precursors,
as found by the restricted search procedure,
are the configurations which are labeled $a$ and $c$ in each figure.
Note that the best found packings of $n=17$ and $n=26$ circles 
shown in the diagrams $b$ and $d$ in Figure~\ref{fig:17a26}  
are obtained from
the best packings found by the restricted search procedure
and 
shown, respectively, in the diagrams $a$ and $c$
in this figure,
using
improvement parameter $\delta = \delta_2$
where  
\begin{equation}
\label{delta2}
\delta_2 =        
2 - 0.5  \sqrt 3  -
 3^{1/4} (2 \sqrt 3 - 1)/(2 \sqrt {4 - \sqrt 3 }).
\end{equation}
Also note that
there are several equivalent ways of the improvement
resulting in the same value of $\delta = \delta_2$.
For example, circle $D$ in 
Figure~\ref{fig:17a26}$b$,
can occupy an alternative position
in which
$D$ contacts
the right side of the rectangle instead of the unlabeled circle
to its left while remaining in contact with circles $B$ and $E$.
This position is also
shown in the figure.
A similar equivalent re-positioning of circle $C$ is shown in
Figure~\ref{fig:17a26}$d$. 

The best found packing of $n=38$ circles
shown in Figure~\ref{fig:37a38}$b$ is obtained from
the best packing found by the restricted search procedure
and 
shown in Figure~\ref{fig:37a38}$a$
using
improvement parameter $\delta = \delta_3$.
The best found packing of $n=58$ circles
shown in Figure~\ref{fig:58a62}$b$ is obtained from
the best packing found by the restricted search procedure
and 
shown in Figure~\ref{fig:58a62}$a$
using
improvement parameter $\delta = \delta_4$.
The values of $\delta_i$ are given in Table~\ref{tab:delta}.
In the packing $a$ shown in Figure~\ref{fig:58a62} and in all the previously
discussed packing diagrams, 
the distances in 
circle-circle or circle-boundary pairs are zero 
whenever the circles in each pair are
 apparently in contact with each other.
The packing $b$ in 
Figure~\ref{fig:58a62} gives an exception to this rule:
the distance between the circles $C$ and $J$ 
and that between the circles $G$ and $K$ is
0.05323824.. of the common circle radius,
perhaps too small to be discerned as positive from the diagram.

Note that the $\delta$-improvement of the configuration
sometimes releases certain circles from the contacts with their
neighbors. The released circles become the so-called {\em rattlers}.
In Figure~\ref{fig:17a26}, circle $A$ becomes a rattler
during the $\delta_2$-conversion of the configuration $c$
into
the configuration $d$.
In Figure~\ref{fig:37a38}, circle $D$ becomes a rattler
during the $\delta_1$-conversion of the configuration $c$
into
the configuration $d$.
Rattlers are represented by
unshaded circles in the packing diagrams.

\begin{table}
\begin{center}
\fbox{
\begin{tabular}{r|c|c} 
$i$&$\delta_i$& defined or used in \\ \hline
   1  & 0.13879028..
           &  Equation~\eqref{delta1}, 
            Figures~\ref{fig:7}$d$, \ref{fig:37a38}$d$   \\
   2  & 0.05728065..   
           &  Equation~\eqref{delta2}, Figures~\ref{fig:17a26}$b$, 
           \ref{fig:17a26}$d$, \ref{fig:57}$b$        \\
   3  & 0.01935364..   
           &  Figures~\ref{fig:37a38}$b$, \ref{fig:31}$b$, \ref{fig:43}$b$,
              \ref{fig:58a62}$d$
  \\
   4  & 0.00403953..   & Figure~\ref{fig:58a62}$b$            
  \\
\end{tabular}
}
\caption{Improvement parameters $\delta_i$}
\label{tab:delta}
\end{center}
\end{table}

In all the examples in Table~\ref{tab:1t62}, the $\delta$-improvement 
of the configuration is localized:
not counting 
the circles possibly used for covering the holes,
all the circles 
involved in the change are located by the boundary on one side.
During the change
the width of the rectangle
decreases by $\delta$ while
the height stays unchanged.
The obtained packings, although they are non-regular, 
are close to their regular-with-holes precursors.
We will call such non-regular packings {\em semi-regular}.
All the non-regular packings listed by their
regular-with-holes representations in 
Table~\ref{tab:1t62} are semi-regular.

We skipped several values of $n$ in 
Table~\ref{tab:1t62}. 
The best found packings obtained for the skipped $n$
show more irregularity than the 
semi-regular packings do.
These excluded from
Table~\ref{tab:1t62}
packings will be called {\em irregular}.
An irregular packing is defined by negation:
it is a packing that cannot be generated
by the described above simple adjustment
where the circles that move are
limited to those covering the holes and
to those located at a side column of a regular-with-holes packing.

We conclude this section with the following observation.
If among the best packings in $R_n$ delivered by the
restricted search there is at least one with holes
or if a non-regular packing of $n$ circles is known
with a smaller perimeter than of those best in $R_n$,
then the packing of $n$ circles
in a rectangle of the minimum perimeter
{\em provably} cannot have a regular pattern.
That is, it cannot be purely hexagonal or purely square-grid
or a hybrid pattern. 
Thus, the optimum packings for each
star-marked semi-regular $n$ in
Table~\ref{tab:1t62} cannot possibly be regular.
We will see in the following section that
the optimum packings for the irregular $n$,
those skipped in
Table~\ref{tab:1t62}, cannot be regular either:
for each skipped $n$ we will present a packing which is better
than the record best in $R_n$.

\section{Results: irregular optimal packings}\label{sec:irreg}
\begin{figure}
\centering
\includegraphics*[width=5.8in]{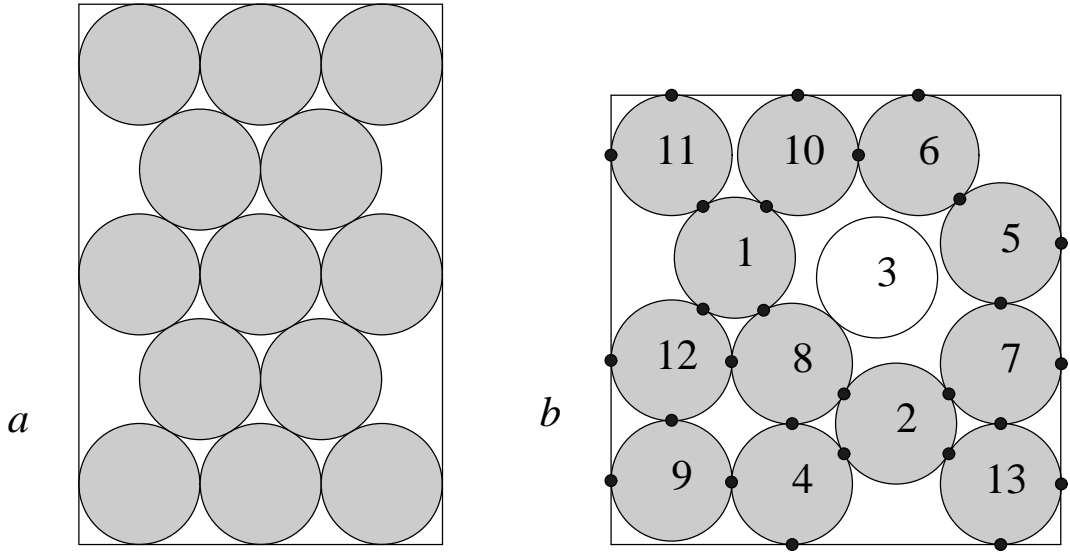}
\caption{Packings of 13 equal circles in rectangles:
$a$) with the smallest perimeter of the enclosing rectangle
among the set $R_{13}$, 
$b$) with the smallest perimeter of the enclosing rectangle
we could find 
}
\label{fig:13}
\end{figure}

\hspace*{\parindent} 
The smallest skipped entry in
Table~\ref{tab:1t62} is $n = 13$.
Figure~\ref{fig:13}$a$ presents the only existing best in $R_{13}$ packing
as found by the restricted search procedure.
The packing has $w = 3$, $h = 5$, $h_- = 2$.
Its perimeter is $P/r = 16 + 8 \sqrt{3} = 29.856406...$.
Since there is no hole in the packing, 
the case $n = 13$ would have qualified as a regular one
and would have been listed as such 
with its parameters $w$, $h$, and $h_-$
in Table~\ref{tab:1t62} 
were it not for the ``compactor'' simulation.
Unexpectedly for us, the ``compactor'' produced 
a better packing!
That packing with the perimeter
$P^{opt}/r =  29.851847510...$ which is smaller
than the perimeter of the packing $a$ 
in Figure~\ref{fig:13} by at least 0.004  
is the packing $b$ shown in the same figure.
The pattern of the packing $b$
in Figure~\ref{fig:13} is truly irregular and non-obvious,
unlike the straight-forward pattern of the packing $a$.
Even the existence of the packing $b$ should not be taken for granted.
By contrast the existence of the packing $a$ 
in Figure~\ref{fig:13} 
can be easily proven by construction
and so can the existence of all the other
regular and semi-regular packings discussed above.

The small black dots in the packing diagram $b$ 
indicate the so-called {\em bonds} or contact points
in circle-circle or circle-boundary
pairs.
A bond
indicates the distance being exactly zero between the pair,
while the absence of a bond in a spot of an apparent contact
indicates the distance being positive, i.e. no contact.
For example, there is no contact between circle 9 and the bottom
boundary in the packing $b$.
(The 13 circles are arbitrarily assigned distinct labels 1 to 13 
in Figure~\ref{fig:13}$b$ to facilitate their
referencing.)
No bond indication is needed in the diagram of the packing $a$
in Figure~\ref{fig:13}
nor in any other regular 
packing diagram
because the points of apparent contacts 
are always the true
contacts in such packings.
In semi-regular packings such non-contacts
do occur, for example 
the one between circles $C$ and $J$
in 
Figure~\ref{fig:58a62}$b$.

With the explicit indication of the bonds,
it is {\em provably} possible to construct the packing 
in Figure~\ref{fig:13}$b$ and this construction is {\em provably} unique,
so the positions of all circles, except the rattler, 
and the rectangle dimensions are
uniquely defined.
The computed horizontal width and vertical heights of the packing
in units equal the circle radius $r$ are
\\
$~~~~~~~~~~~~~~~~~~~~~~~W/r =  5.463267269314... ~~~~H/r =  5.462656485780...$
\\
which implies the perimeter value given above
and 
$L/S = 1.000111810716...$
so the rectangle is almost a square to within about 0.01\% 

Note that 
when in \cite{LG} we minimized the area of the rectangle
we were unable to find packings better than,
in the present paper terminology,
either regular or semi-regular.
The case $n = 13$ is a violation
of such structure for the case of 
minimizing the perimeter.
Suspicious of other such violations,
we ran
many more tries 
of the ``compactor'' simulation for $n=14,15,16,17,18,19$ and 20.
No violation was detected.
All these cases seem to be either regular or semi-regular
as presented in Table~\ref{tab:1t62}.
But for $n = 21$ we encountered another gross 
violation of regularity.

\begin{figure}
\centering
\includegraphics*[width=5.8in]{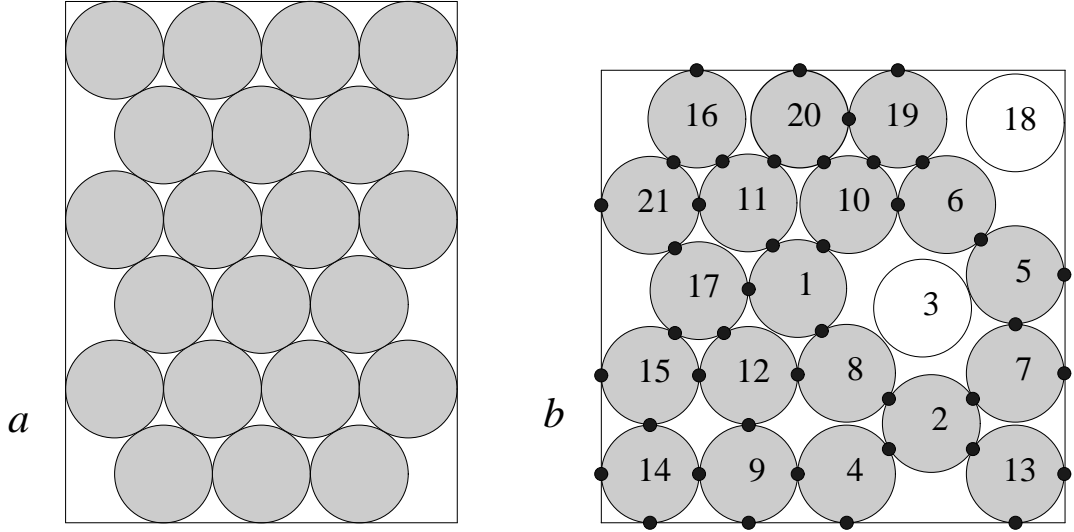}
\caption{Packings of 21 equal circles in rectangles:
$a$) with the smallest perimeter of the enclosing rectangle
among the set $R_{21}$, 
$b$) with the smallest perimeter of the enclosing rectangle
we could find 
}
\label{fig:21}
\end{figure}
The case of $n = 21$ is similar to that of $n = 13$.
Here again, the restricted search procedure
delivers the best in $R_{21}$ packing (shown in Figure~\ref{fig:21}$a$)
with $w = 4$, $h = 6$, $h_- = 3$ and 
perimeter $P/r = 20 + 10 \sqrt{3} = 37.3205081..$.
The packing is regular and looks very different from the best packing
found by the simulation (shown in Figure~\ref{fig:21}$b$),
the latter
with the perimeter $P^{opt}/r = 37.309294229..$ which is smaller
than the perimeter of the packing 
shown in Figure~\ref{fig:21}$a$ by at least 0.01.  
Same as in the case of $n = 13$, 
the best found packing for $n = 21$ 
exhibits a rather irregular structure,
which makes 
the existence of the packing non-obvious.
With the bonds shown in Figure~\ref{fig:21}$a$,
it is possible to {\em prove} the existence
of the packing by construction and 
it is possible to {\em provably} uniquely determine the position
of the circles, except the rattlers, 
the width, height and the $L/S$ ratio of the rectangle:
\\
$~~~~~~~~~~~~~W/r =  7.433745175630.. ~~~H/r =  7.220901938764.. ~~~ L/S = 1.029475990489..$

The patterns of both irregular best found packings,
while being dissimilar to all the other
conjectured optimum packings considered thus far,
show some similarity between themselves.
This similarity
is emphasized by the circle labeling.
Labels 1 to 13 
in Figure~\ref{fig:21}$b$ are assigned
to the circles that occupy the positions in
that figure which are similar to the corresponding
circles 1 to 13 in Figure~\ref{fig:13}$b$.
The similarity, however, is not perfect.
For example, 
the bond between circle 9 and the bottom boundary
in Figure~\ref{fig:21}$b$ does not find its counterpart 
in Figure~\ref{fig:13}$b$.

\begin{figure}
\centering
\includegraphics*[width=6.0in]{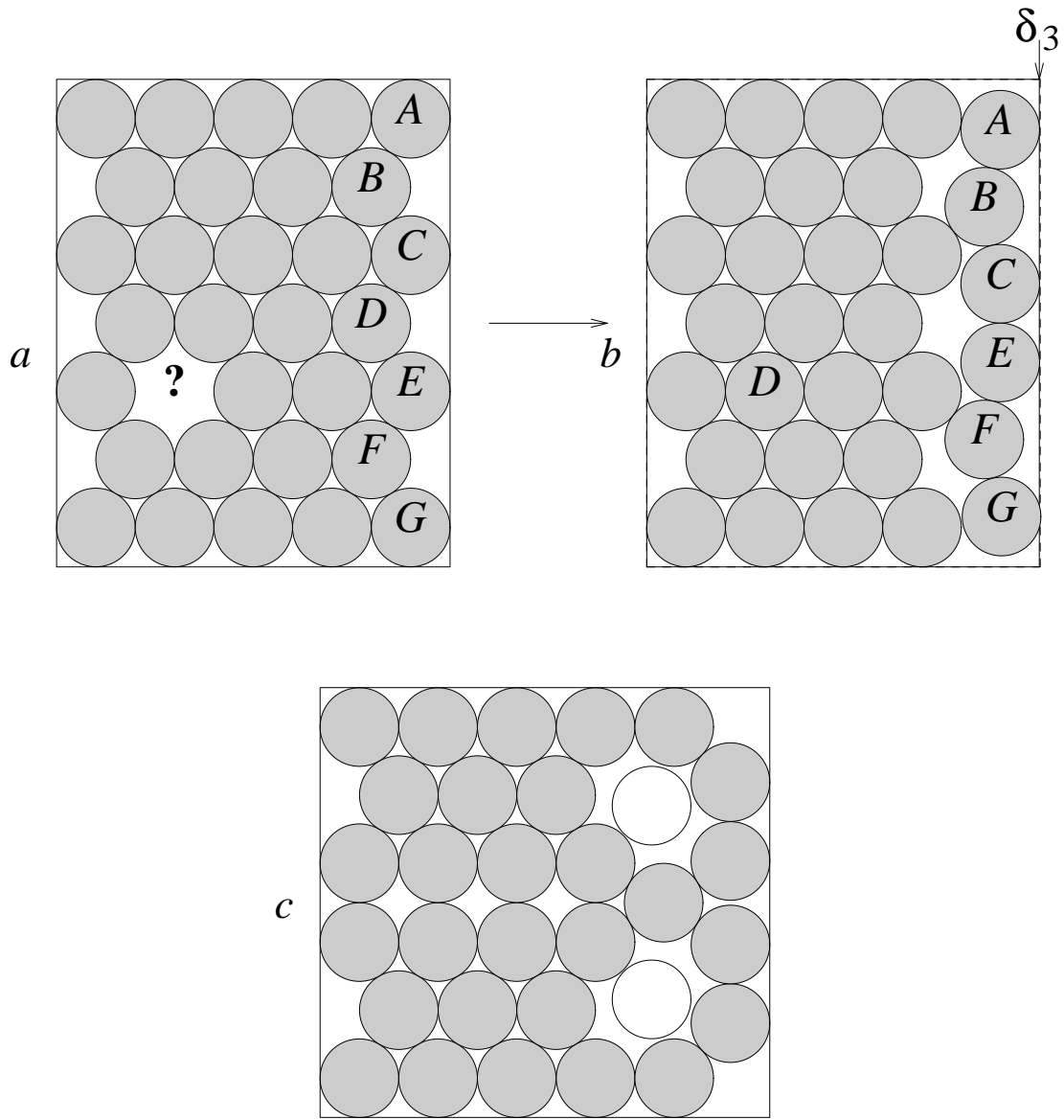}
\caption{Packings of 31 equal circles in rectangles:
$a$) with the smallest perimeter of the enclosing rectangle
among the set $R_{31}$, 
$b$) $\delta_3$-improved packing $a$,
$c$) with the smallest perimeter of the enclosing rectangle
we could find 
}
\label{fig:31}
\end{figure}

\begin{figure}
\centering
\includegraphics*[width=6.1in]{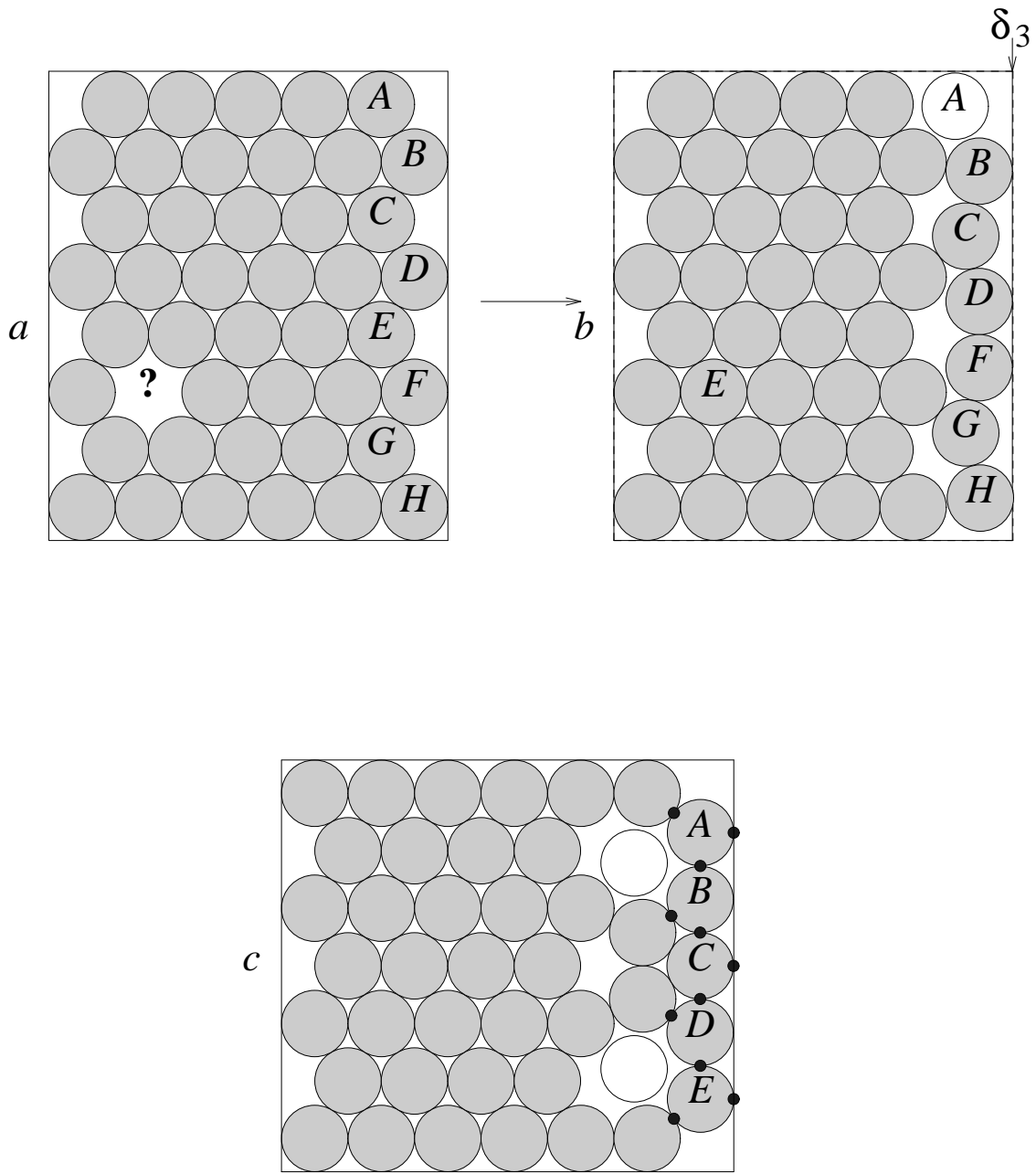}
\caption{Packings of 43 equal circles in rectangles:
$a$) with the smallest perimeter of the enclosing rectangle
among the set $R_{43}$, 
$b$) $\delta_3$-improved packing $a$,
$c$) with the smallest perimeter of the enclosing rectangle
we could find. The black dots indicate bonds
of the labeled circles in the packing $c$
}
\label{fig:43}
\end{figure}

\begin{figure}
\centering
\includegraphics*[width=6.3in]{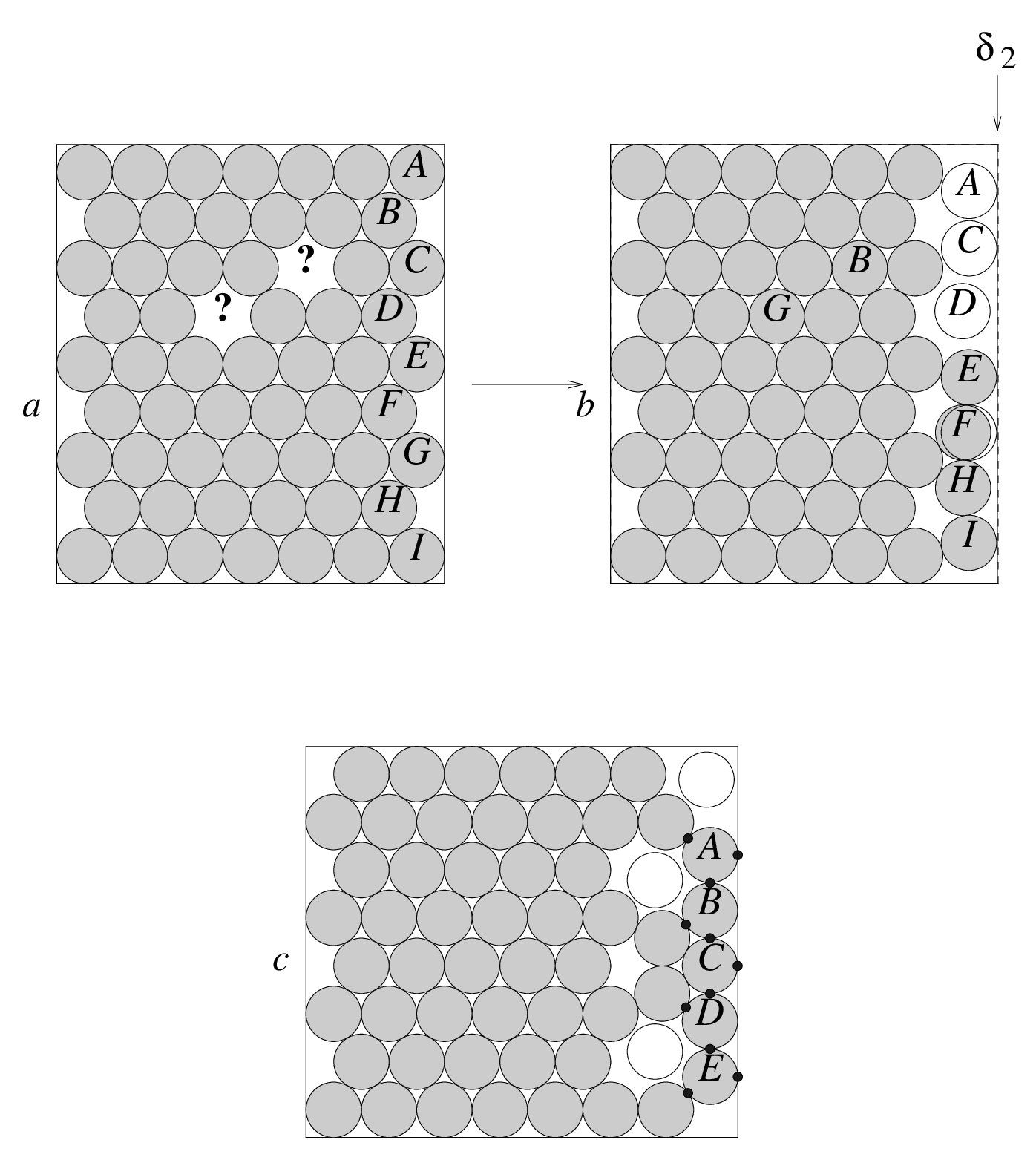}
\caption{Packings of 57 equal circles in rectangles:
$a$) with the smallest perimeter of the enclosing rectangle
among the set $R_{57}$, 
$b$) $\delta_2$-improved packing $a$,
$c$) with the smallest perimeter of the enclosing rectangle
we could find. 
Alternative position of circle $F$ is shown in packing $b$.
The black dots indicate bonds
of the labeled circles in the packing $c$
}
\label{fig:57}
\end{figure}

The remaining skipped entries in Table~\ref{tab:1t62}
are $n = 31$, 43, and 57. 
For these three values of $n$, unlike $n=13$ or 21,
the best in $R_n$ packings all have holes,
as seen in Figures~\ref{fig:31}$a$,
\ref{fig:43}$a$, and
\ref{fig:57}$a$, 
and hence avail themselves for $\delta$-improvements.

The improved packings
labeled $b$ in these three figures
have perimeters, respectively
\\
$P/r = 
12(2+\sqrt{3}) - 2 \delta_3
=
44.74590240843..
$ for 31 circles,
\\
$P/r = 
14(2+\sqrt{3}) - 2 \delta_3
=
52.21000402357..
$ for 43 circles,
\\
$P/r = 
16(2+\sqrt{3}) - 2 \delta_2
=
59.59825161939..
$ for 57 circles.

Those improved perimeters still exceed the 
perimeters of the corresponding best packings found, 
which happen to be irregular, namely
\\
for 31 circles $P^{opt}/r = 
44.7095500424198..
$ is smaller than $P/r$ by at least 0.035,
\\
for 43 circles $P^{opt}/r = 
51.99029827020367..
$ is smaller than $P/r$ by at least 0.2,
\\
for 57 circles $P^{opt}/r = 
59.4543998853414..
$ is smaller than $P/r$ by at least 0.14.

The patterns of the irregular packings $c$ in 
Figures~\ref{fig:31},
\ref{fig:43}, and
\ref{fig:57} somewhat resemble each other, especially
the packings of 43 and 57 circles.
Moreover, the best found packing of 43 circles in 
Figure~\ref{fig:43}$c$ is an exact subset in
the best found packing of 57 circles in
Figure~\ref{fig:57}$c$.
In both packings, the dots attached to 5 circles labeled $A$ to $E$
indicate the bonds of these circles.
Thus, for example, the circles $B$ and $D$ do contact the circle $C$,
but do not contact the right side of the rectangle, where
the gap is 
0.00957...
of the circle radius, too small to be 
discerned in the diagrams.
Similarly, the circle $C$ does contact the right side of the rectangle
but does not contact either of the two unlabeled
circles immediately at $C$'s left. 
Between the pairs that do not include at least one labeled circle
the bonds exist in the obvious places,
and they are not specifically indicated by dots
in the figures.
As before, it is possible to {\em prove} the existence of
the irregular packings $c$ in
Figures~\ref{fig:31},
\ref{fig:43}, and
\ref{fig:57}. 
This existence might be non-obvious, especially,
for the latter two packings.

\section{Double optimality and related properties}\label{sec:dopt}
\begin{figure}
\centering
\includegraphics*[width=6.0in]{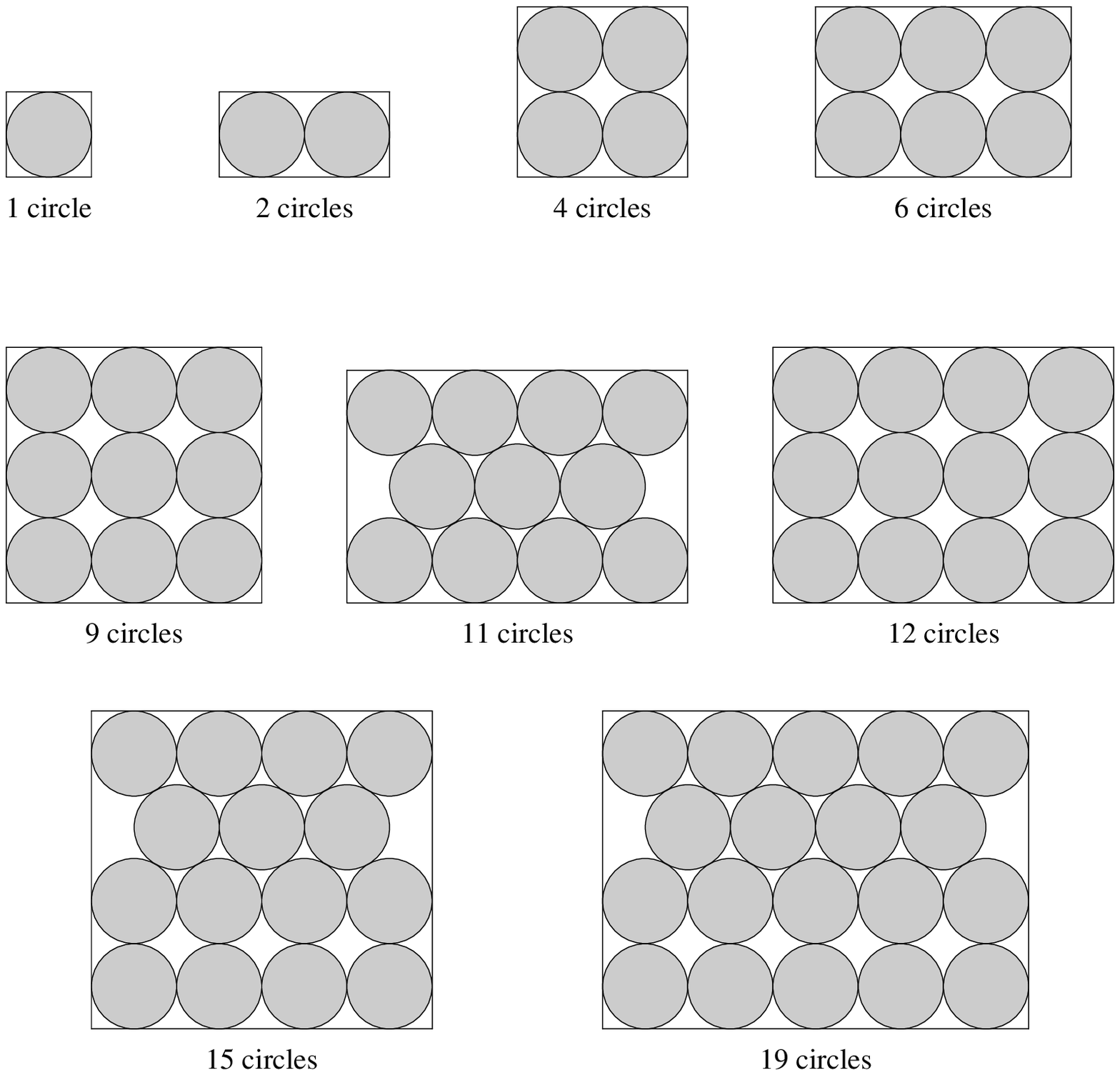}
\caption{The rectangle that encloses
each of these packings has the smallest perimeter and the smallest area
that we could find for their number of circles
}
\label{fig:same}
\end{figure}

\hspace*{\parindent} 
As mentioned in Introduction,
minimizing the perimeter and minimizing the area
of the rectangle
lead to 
generally different optimal packings
for the same number of circles $n$,
if the rectangle aspect ratio is variable.
Are there packings
optimal under both criteria at the same time?
Figure~\ref{fig:same} displays 
such conjectured double optimal packings 
that were obtained by comparing the list of smallest area
packings reported in \cite{LG} with that of
the smallest perimeter packings reported here.
For a particular $n$ and a particular optimality criterion, 
there may be several equivalent optimum packings.
For example, 
two minimum-perimeter packings exist for $n=11$ according
to Table~\ref{tab:1t62} and two minimum-area packings exist
for $n=15$ according to \cite{LG}.
However, no more than one packing was found to be double optimal
for any $n$.

\begin{figure}
\centering
\includegraphics*[width=6.1in]{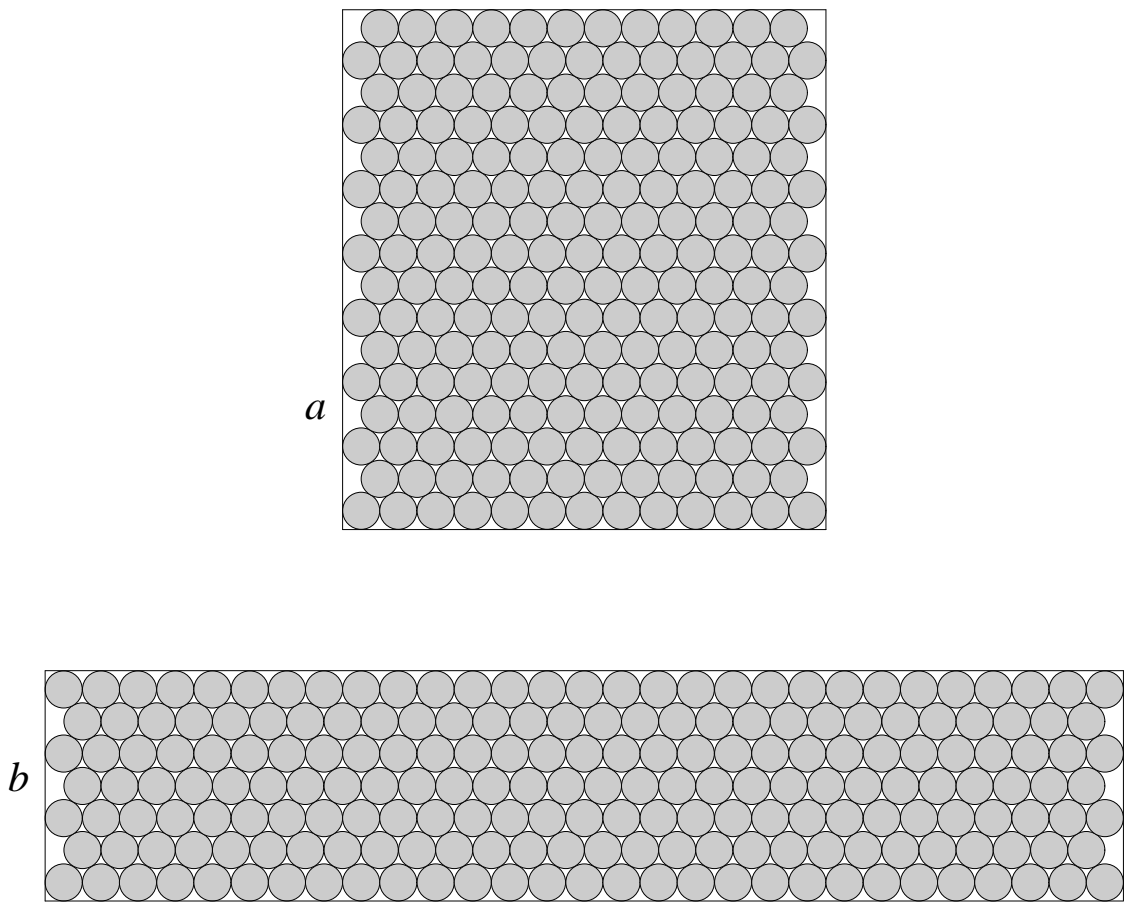}
\caption{The best packings found for $n=200$ equal circles
in rectangles of a variable shape: $a$) under the criterion
of the minimum perimeter, $b$) under the criterion
of the minimum area
}
\label{fig:200}
\end{figure}

Since we will show
in Section~\ref{sec:prf}
that for $n \rightarrow \infty$ the $L/S$ ratio
of the minimum-perimeter rectangles tends to 1,
and it is conjectured in this case 
that the $L/S$ ratio for the minimum-area rectangles tends to
$2 + \sqrt{3}$ (see \cite{LG}),
then, conjecturally,
there may be only a finite number of double optimal packings.
In fact, we believe  Figure~\ref{fig:same} lists 
all double optimal packings.
For larger $n$, best rectangular shapes found
under the two optimality criteria become 
noticeably different from each other,
e.g., see the best packings found under either criteria
for $n = 200$ in Figure~\ref{fig:200}.

Consider a configuration $C^{opt}$ of $n$ equal circles
which supplies the
global minimum perimeter for its enclosing rectangle
and suppose the $C^{opt}$
happens to be {\em not} the one that supplies
the global minimum to the area of the enclosing rectangle.
We believe, however, that the rectangle of this $C^{opt}$ still 
holds a record of being a rectangle of the minimum area,
albeit locally.
Specifically, if we vary slightly the ratio $L/S$ of the rectangle
around its value given by $C^{opt}$ and for each such $L/S$
find the rectangle of the densest possible packing of $n$ unit-radius circles
(this rectangle possesses
both the minimum perimeter and the minimum area for its value $L/S$),
then the area of the rectangle for the configuration $C^{opt}$ 
will turn out to be the minimum among the areas of all those varied rectangles.
We also believe that the statement which is obtained from
the statement above by the interchange of
the minimum-perimeter criterion with the minimum-area criterion is also
true. 
That is, 
a configuration that delivers the 
global minimum of the area of the enclosing rectangle
also holds a record of supplying the minimum of its perimeter, 
though perhaps only locally.
Can the two sets of configurations,
those that deliver the local minima for the rectangle area 
and those that deliver the local minima for the rectangle perimeter,
be the same sets?

We have tested the former statement 
(that the global minimum-perimeter optimality implies the local
minimum-area optimality) numerically
for some values of $n$. 
If this conjecture is true,
then the minimum-perimeter packings for some values of $n$ have to be
of a higher density than the densest packings of $n$ equal circles
in a square.
Cases $n=13$ and $n=21$ seem to be
such occurrences.
The configurations of 13 and 21 circles with
the smallest found rectangular perimeter shown
in Figures~\ref{fig:13}$b$ and \ref{fig:21}$b$ resemble the
respective best found packings of 13 and 21 equal circles in a square,
see for example, \cite{GL2}.
The only visible difference between the pairs for each $n$ 
is in the positions of the bonds. 
The $L/S$ ratios in either minimum-perimeter packing
is very close to 1, so
each respective densest packing in a square is a local neighbor
of the corresponding minimum-perimeter packing.
It can be verified that each of the two densest-in-a-square packings
in \cite{GL2}
has a lower density than that of
its minimum-perimeter counterpart reported here.

\section{Minimum-perimeter packings for larger $n$}\label{sec:largn}
\hspace*{\parindent} 
Because the results of the restricted search, available for all
$n \le 5000$, were
supported by simulation only for all $n \le 62$,
our conjectures become more speculative for larger $n$.

For $n \le 62$
the restricted search reliably predicts the best packings
found by the simulation
except those of the form 
\begin{equation}
\label{k(k+1)+1}
n = k(k+1)+1,
\end{equation}
where $3 \le k \le 7$.
All those predicted packings
happen to be either regular or semi-regular.
The optimal packings for the values of $n$ of the form
\eqref{k(k+1)+1}
appear to be exceptionally irregular.
For which $n > 62$ are the best packings similarly irregular?

Our previous experience for packing equal circles
in various shapes suggests that
the numbers of circles which result in
exceptionally ``bad'' or irregular optimal packings
often follow immediately after the numbers which result 
in exceptionally ``good'' or regular optimal packings.
For example, 
triangular numbers of circles $n=k(k+1)/2$
arrange themselves optimally in regular triangular patterns
inside equilateral triangles and
in \cite{GL1} we observed that
the optimal arrangements 
of $n=k(k+1)/2+1$ circles
look irregular and
disturbed.

The minimum-perimeter packings of $n = k(k+1)$ equal circles
inside rectangles with a variable aspect ratio, as 
conjectured 
by the restricted search for 
\begin{equation}
\label{333}
4 \le k \le 33 ,
\end{equation}
are
regular hexagonal arrangements of 
$h = k+1$ alternating rows 
with $w = k$ circles in each row.
For $k > 33~(n > 1122)$ this
regular pattern with $w=k$, $h=k+1$ does not serve as the optimum.
The latter statement is proven at least for $n \le 5000$.
Thus we speculate that the 30 values of $k$ in \eqref{333} also
might correspond to the 30 cases 
of $n$ as computed by formula \eqref{k(k+1)+1}
in each of which the minimum-perimeter
packing of $n$ circles in rectangles with variable aspect ratio
is irregular.

The irregular minimum-perimeter packings,
probably, also occur
for some $n$ which are
{\em not}
of the form \eqref{k(k+1)+1}.
We believe $n = 66$ is the smallest such $n$.
In fact, $n = 66$ is the smallest one with the properties:
(A) it is not of the form \eqref{k(k+1)+1} for any integer $k$,
(B) the best in $R_n$ packing, as delivered by the restricted search,
has 
$h = 9$ alternating rows,
$h_{-} = 4$ of which are one circle shorter, and 
$v = 2$ holes.
The smallest $n$ for which (B) holds is
$n = 57$.
The best in $R_{57}$ packing is shown
in Figure~\ref{fig:57}$a$.
By attaching
a column of 9 alternating circles at the left of
Figure~\ref{fig:57}$a$ we obtain a diagram of the 
best in $R_{66}$ packing.
The best in $R_{57}$ packing
can be $\delta_2$-improved as shown
in Figure~\ref{fig:57}$b$.
The best in $R_{66}$ packing
can be $\delta_2$-improved in the same way.
In either one of these $\delta_2$-improved packings
all the circles can be ``unjammed''
so that there would be no contacts among them or with the boundary.
Hence the perimeter of either one can be further
reduced by subsequent ``compaction'' of the rectangle.

Our ``compactor'' simulation suggests
that for $n = 57$ 
the irregular packing
thus obtained
does not overtake the conjectured
optimum packing shown in
Figure~\ref{fig:57}$c$.
However, we believe
that its analogue for $n = 66$ might just be the optimum packing
and it would be irregular.
Approximately,
each of the two ``compacted'' irregular packings might look similar
to the packing in Figure~\ref{fig:57}$b$.
Unfortunately,
obtaining their exact patterns,
including the identification of the bonds,
proved to be beyond our
current computing capabilities.
(The hardness of the computations might indicate existence
of several local minima
near the ``unjammed'' configurations.)
Note that the pattern of the least perimeter
packing for $n = 66$
would probably differ substantially 
from the patterns of the irregular least perimeter packings
for $n = 13, 21, 31, 43$, and 57,
those conjectured
in Section~\ref{sec:irreg}.

The chance to encounter a non-regular $n$ (which 
by definition
has to correspond to either a semi-regular or irregular optimal packing)
increases quickly with $n$.
In Table~\ref{tab:5000} we list the 
packings found by the restricted search procedure for several segments
of consecutive $n$.
The segments are arbitrarily selected within the set $62 < n \le 5000$.
%
\begin{table}
\begin{center}
\fbox{
\begin{tabular}{r|r|r|r|r|r||r|r|r|r|r|r||r|r|r|r|r|r} 
$n$&$w$&$h$&$h_{-}$&$s$&$v$&
$n$&$w$&$h$&$h_{-}$&$s$&$v$&
$n$&$w$&$h$&$h_{-}$&$s$&$v$ \\ \hline
 101  &  9 & 12 &  6 &  0 & 1         &  501 & 21 & 24 &  0 &  0 &  3
                                    & 2001 & 44 & 46 & 23 &  0 &  0     \\
 102  &  9 & 12 &  6 &  0 & 0         &  502 & 21 & 24 &  0 &  0 &  2
                                    & 2002 & 41 & 49 &  0 &  0 &  7     \\
 103  & 10 & 11 &  5 &  0 & 2         &  503 & 21 & 24 &  0 &  0 &  1
                                    & 2003 & 41 & 49 &  0 &  0 &  6     \\
 104  & 10 & 11 &  5 &  0 & 1         &  504 & 21 & 24 &  0 &  0 &  0
                                    & 2004 & 41 & 49 &  0 &  0 &  5     \\
 105  & 10 & 11 &  5 &  0 & 0         &  505 & 22 & 23 &  0 &  0 &  1
                                    & 2005 & 41 & 49 &  0 &  0 &  4     \\
 106  & 10 &  9 &  4 &  2 & 0         &  506 & 22 & 23 &  0 &  0 &  0
                                    & 2006 & 41 & 49 &  0 &  0 &  3     \\
      & 11 &  9 &  4 &  1 & 0         & *507 & 20 & 26 & 13 & 0 & 0
                                    & 2007 & 41 & 49 &  0 &  0 &  2     \\
 107  &  9 & 12 &  0 &  0 & 1         &  508 & 21 & 25 & 12 & 0 & 5
                                    & 2008 & 41 & 49 &  0 &  0 &  1     \\
 108  &  9 & 12 &  0 &  0 & 0         &  509 & 21 & 25 & 12 & 0 & 4
                                    & 2009 & 41 & 49 &  0 &  0 &  0     \\
 109  & 10 & 11 &  0 &  0 & 1         &  510 & 21 & 25 & 12 & 0 & 3
                                    & 2010 & 42 & 48 &  0 &  0 &  6     \\
 110  & 10 & 11 &  0 &  0 & 0         &  511 & 21 & 25 & 12 & 0 & 2
                                    & 2011 & 42 & 48 &  0 &  0 &  5     \\
 ...  & .. & .. &  . &  . & .         &  ... & .  & . & . & . & .
                                    &  .. & ..& ..& ..& ..& ..          \\
 251  & 14 & 18 &  0 &  0 & 1         & 1001 & 30 & 34 & 17 & 0 & 2
                                    & 4991 & 64 & 78 &  0 &  0 &  1     \\
 252  & 14 & 18 &  0 &  0 & 0         & 1002 & 30 & 34 & 17 & 0 & 1
                                    & 4992 & 64 & 78 &  0 &  0 &  0     \\
 253  & 15 & 17 &  0 &  0 & 2         & 1003 & 30 & 34 & 17 & 0 & 0
                                    & 4993 & 68 & 74 & 37 &  0  & 2     \\
 254  & 15 & 17 &  0 &  0 & 1         & 1004 & 31 & 33 & 16 & 0 & 3
                                    & 4994 & 68 & 74 & 37 &  0  & 1     \\
 255  & 15 & 17 &  0 &  0 & 0         & 1005 & 31 & 33 & 16 & 0 & 2
                                    & 4995 & 68 & 74 & 37 &  0  & 0     \\
 256  & 16 & 16 &  0 &  0 & 0         & 1006 & 31 & 33 & 16 & 0 & 1
                                    & 4996 & 65 & 77 &  0 &  0  & 9     \\
 257  & 14 & 19 &  9 &  0 & 0         & 1007 & 31 & 33 & 16 & 0 & 0
                                    & 4997 & 65 & 77 &  0 &  0  & 8     \\
 258  & 15 & 18 &  9 &  0 & 3         & 1008 & 28 & 36 &  0 & 0 & 0
                                    & 4998 & 65 & 77 &  0 &  0  & 7     \\
 259  & 15 & 18 &  9 &  0 & 2         & 1009 & 29 & 35 &  0 & 0 & 6
                                    & 4999 & 65 & 77 &  0 &  0  & 6     \\
 260  & 15 & 18 &  9 &  0 & 1         & 1010 & 29 & 35 &  0 & 0 & 5
                                    & 5000 & 65 & 77 &  0 &  0  & 5     \\
 ...  & .. & .. &  ..&  ..& ..        & .... & .. & .. &  ..& ..& ..
                                    & ... & . & . & . & . &  . \\
\end{tabular}
}
\caption{Packings of $n$ circles in rectangles 
of the smallest perimeter as found by the restricted search
for several contiguous segments of $n$ that were
arbitrarily selected
within the set $62 < n \le 5000$.
The packings shown are either regular (with $v = 0$)
and then they 
all
are believed to be globally optimal, with the exception of the case $n=507$,
or they are regular with holes
(with $v > 0$) and then they can be $\delta$-improved into
semi-regular packings.
The exceptional entry $n = 507 = 22\times23+1$ is marked with a star.
The optimal packing for 507 circles might 
be of an unknown irregular pattern
}
\label{tab:5000}
\end{center}
\end{table}
The structure of Table~\ref{tab:5000} is similar to that
of Table~\ref{tab:1t62}, except that 
an additional column is provided for
the number of holes $v$.
The entries $n$ with $v > 0$ are frequent
in Table~\ref{tab:5000}, unlike
Table~\ref{tab:1t62}.
Also, we choose to skip the $\delta$ column here.
(Determining and presenting the higher-order $\delta_i$
would involve many details
exceeding the reasonable limits for this paper.)

The discussion above suggests
that, perhaps, 
some of the entries with 
multiple holes, $v > 1$,
correspond to irregular packings,
if their $\delta$-improvements
avail themselves 
to further improvements same as the packing in 
Figure~\ref{fig:57}$b$.\footnote{
As reported in \cite{LG},
for some $n \ge 393$ the configurations with the least area
among the set $R_n$ have multiple holes.
In particular,
the smallest $n$ for which 
property (B) holds for
the least rectangular area configuration
among the set $R_n$ 
is $n = 453$. 
The non-zero parameters of the configuration 
are $w~=~51, h~=~9, h_{-}~=~4, v~=~2$.
The packing of $453$ congruent circles which delivers
the global minimum to the area of the enclosing rectangle
is, probably, irregular.
}

The larger values of $n$ that would correspond to regular packings 
become rare.
However, we do not believe regular $n$ eventually disappear.
In other words, we do not believe the largest such $n$ exists.
For the infinite sequence of values of $n$,
all those of the form 
\begin{equation}
\label{square1}
n=\frac{1}{2}(a_k+1)(b_k+1),~~ k = 1,2,..
\end{equation}
with
\begin{equation}
\label{square2}
a_1 = 1, a_2 = 3, a_{k+2} = 4a_{k+1} - a_k,~~~
b_1 = 1, b_2 = 5, b_{k+2} = 4b_{k+1} - b_k,~~k = 1,2,...
\end{equation}
the minimum-perimeter
packings are probably regular.
The fractions
$a_k  /b_k$ are (alternate) convergents to
$1/ \sqrt 3$, and it has been conjectured by Nurmela et al.
\cite{NOR} that for these $n$, a ``nearly'' hexagonal packing of
$n$ circles in a square they describe is in fact optimal,
that is, it has the largest possible density among those in a square.
The beginning terms $n = 12, 120, 1512$ 
of sequence \eqref{square1} are within the range $n \le 5000$
for which we exercised the restricted search procedure.
For these three $n$'s
the search delivers regular patterns 
as the minimum-perimeter packing in $R_n$:
a $4 \times 3$ square-grid for $n=12$,
the hexagonal packing
with $w = 10$, $h=12$ for $n=120$ and 
that
with $w = 36$, $h=42$ for $n=1512$.
Note that the latter two are
the {\em exact} hexagonal packings
and they have densities larger than that of the corresponding
``nearly'' hexagonal packings, those that are best in a square
according to the conjecture in \cite{NOR}.
We believe the same relation between
the best packings of $n$ circles in a square and the minimum-perimeter
packings of $n$ circles in rectangles continues
for all larger $n$ of the form \eqref{square1},  \eqref{square2}.

Increasing the value of $n$ seems to diminish 
such phenomena as dimorphism and hybrid packings.
Dimorphism of the optima
is the existence of two different optimal rectangular shapes.
The dimorphism occurs within the interval $1 \le n \le 62$
for $n = 11$, 19, 28, 29, 40 and 53,
see Table~\ref{tab:1t62}.
The hybrid packings, 
identifiable in the table by 
$h$ and $s$ 
being both positive,
occur for $n = 11$, 15, 19, 24, 28, 29, 34, 40, 47, 53, and 61.
Among the $n$'s selected for Table~\ref{tab:5000}, 
the two different optima exist
only for $n = 106$. 
Both optima also happen to be hybrids
and no other hybrid occurs in Table~\ref{tab:5000}.
The remaining cases of dimorphism and/or hybrid packings
among $1 \le n \le 5000$ are all listed 
in Table~\ref{tab:pm}.
\begin{table}
\begin{center}
\fbox{
\begin{tabular}{r|r|r|r|r||r|r|r|r|r||r|r|r|r|r} 
$n$&$w$&$h$&$h_{-}$&$s$&
$n$&$w$&$h$&$h_{-}$&$s$&
$n$&$w$&$h$&$h_{-}$&$s$  \\ \hline
  69  & 8 & 7 & 3 & 2          & 151 &12 &11 & 5 & 2 
                                    & 298 &17 &17 & 8 & 1  \\
      & 9 & 7 & 3 & 1          &     &13 &11 & 5 & 1 
                                    &     &18 &17 & 8 & 0  \\
  78  & 9 & 7 & 3 & 2          & 176 &13 &13 & 6 & 1 
                                    & 316 &18 &17 & 8 & 1  \\
  86  & 9 & 9 & 4 & 1          &     &14 &13 & 6 & 0 
                                    & 371 &19 &19 & 9 & 1  \\
      &10 & 9 & 4 & 0          & 190 &14 &13 & 6 & 1 
                                    &     &20 &19 & 9 & 0  \\
  96  &10 & 9 & 4 & 1          & 233 &15 &15 & 7 & 1
                                    & 452 &21 &21 &10 & 1  \\
 127  &11 &11 & 5 & 1          &     &16 &15 & 7 & 0 
                                    &     &22 &21 &10 & 0  \\
      &12 &11 & 5 & 0          & 249 &16 &15 & 7 & 1 
                                    & 541 &23 &23 &11 & 1  \\
 139  &12 &11 & 5 & 1          &     &   &   &   &   
                                    &     &24 &23 &11 & 0  \\
\end{tabular}
}
\caption{Packings of $n$ circles in a rectangle 
of the smallest found perimeter which are hybrid
and/or exist in two differently shaped rectangles.
All such cases $n > 62$ are listed in this table, except $n = 106$,
which is listed in Table~\ref{tab:5000}
}
\label{tab:pm}
\end{center}
\end{table}
For either phenomenon, dimorphism or hybrid packings,
the largest $n$ for which the phenomenon still occurs
appears to be $n = 541$.
For comparison:
both phenomena also occur for the criterion of the minimum area
in \cite{LG}, but both end much sooner, the largest
$n$ for which either phenomenon still occurs appears to be $n = 31$.

\section{Optimal rectangles are asymptotically square}\label{sec:prf}
\hspace*{\parindent}
In this section we will show that as $n$ goes to infinity, the ratio
of $L/S$ for the minimum perimeter rectangle in which $n$ congruent circles
can be packed tends to 1. This will follow from the following 
considerations. For a compact, convex subset $X$ of the Euclidean
plane, define the {\em packing number} $p(X)$ to be the cardinality of
the largest possible set of points within $X$ such that the distance
between any two of the points is at least 1. 
Let $A(X)$ be the area and $P(X)$ be the perimeter of $X$.
The following result of
Oler (see \cite{Oler}, \cite{FG}) bounds $p(X)$:

{\bf Theorem:}
$$p(X) \le \frac{2}{\sqrt 3} A(X) + \frac{1}{2} P(X) +1.$$
It is easy to prove the following lower bound on the packing number for
a square $S(\alpha)$ of side $ \alpha$:

{\bf Fact}:  $$p(S(\alpha)) \ge \frac{2}{\sqrt 3} \alpha^2$$.

Suppose $R$ is an optimal rectangle with side 
lengths $m + \epsilon$ and $m - \epsilon$.
Thus, $R$ has perimeter $4m$ and area $m^2 - \epsilon^2$. By the preceding
upper bound on $p(R)$ we have

$$p(R) \le \frac{2}{\sqrt 3}(m^2-\epsilon^2) + \frac{1}{2}(4m) + 1.$$

Since $R$ is optimal, then we must have $p(R) \ge p(S(m))$. 
This implies that

$$\frac{2}{\sqrt 3}(m^2-\epsilon^2) + 2m + 1 \ge \frac{2}{\sqrt 3}m^2.$$

From this it follows that 

$$\epsilon \le ( \frac{\sqrt 3}{2}(2m+1) )^{1/2},$$
which is of a lower order than $m$, the order of the side lengths.

It is now straight-forward to convert this inequality to one for packing
circles rather than points, and our claim is proved. As an example, if
$m = 1000$ then $\epsilon < 41.623$..., so that the ratio $L/S \le 1.0869.$


\end{document}